# Reordering of the Logistic Map with a Nonlinear Growth Rate

Dominique Delcourt

*Laboratoire de Physique des Plasmas*
CNRS-UPMC, Paris, France

**Abstract.** In the well known logistic map, the parameter of interest is weighted by a coefficient that decreases linearly when this parameter increases. Since such a linear decrease forms a specific case, we consider the more general case where this coefficient decreases nonlinearly as in a hyperbolic tangent relaxation of a system toward equilibrium. We show that, in this latter case, the asymptotic values obtained via iteration of the logistic map are considerably modified. We demonstrate that both the steepness of the nonlinear decrease as well as its upper and lower boundaries significantly alter the bifurcation diagram. New period doubling features and transitions to chaos appear, possibly leading to regimes with small periods. Computations with a variety of parameter values show that the logistic map can be significantly reordered in the case of a nonlinear growth rate.

# 1. Introduction

When one is to model the transition of some physical system from one configuration to another, the functional form that is adopted to describe this transition is of great importance. One may write this transition as : $\mathscr{E}(t) = \mathscr{E}_{initial} + f(t)[\mathscr{E}_{final} - \mathscr{E}_{initial}]$ where $\mathscr{E}(t)$ is the configuration of the system at a given time $t$, $\mathscr{E}_{initial}$ and $\mathscr{E}_{final}$ are the initial and final configurations, respectively, and $f(t)$ is a functional form that varies between 0 at the beginning of the transition and 1 at the end of the transition. If one chooses a form $f(t)$ that varies linearly with time, this leads to specific, and possibly unrealistic, features. As an example, in plasma physics, if one studies the relaxation of magnetic field lines as a result of some instability from a given geometry to a less stretched one (a process referred to as "dipolarization") and if one adopts the above description for the magnetic field vector **B** (viz., **B**$(t)$ = **B**$_{initial}$ + $f(t)\times$[**B**$_{final}$ - **B**$_{initial}$]), a linear functional form $f(t)$ leads to a constant induced electric field according to Maxwell's equation : $\nabla \times \mathbf{E} = -\partial \mathbf{B}/\partial t$. Such a constant induced electric field abruptly imposed at the beginning of the magnetic transition has no physical meaning and actually leads to artifacts in numerical simulations of the plasma response such as abrupt jumps of the magnetic moment of charged particles. This is not the case if one uses a nonlinear functional form $f(t)$ such as a polynomial of degree 5 that features zero time derivatives at the beginning and at the end of the magnetic transition [e.g., *Delcourt et al.*, 1990].

In the present study, we examine the effect of such a nonlinear functional form $f(t)$ on the logistic map. In the seminal work of *May* [1976], this logistic map has a linear growth rate that regularly decreases as the parameter of interest increases. In this study, we consider a growth rate that changes more slowly than the linear one near minimum and maximum values of the parameter of interest, and more rapidly than the linear one at intermediate values as is the case in the relaxation of a system toward equilibrium. We will show that such a nonlinear growth rate has a significant impact on the long term behavior and thus on the bifurcation diagram. In section 2, we first describe the recurrence relation of interest in the present study. Numerical results are presented in section 3 and put into perspective in section 4.

## 2. Recurrence Relation

Since the discovery of its unique properties in the seminal work of *May* [1976], the logistic map has been a widely investigated topic. This equation writes :

$$Y_{n+1} = C \times (1-Y_n) \times Y_n \qquad (1)$$

where C is a constant (also referred to as the bifurcation parameter). In the above equation, this constant C is weighted by a normalized growth rate (denoted by γ hereinafter) that varies linearly with $Y_n$, viz., $\gamma(Y_n) = 1-Y_n$. This growth rate is such that $\gamma(Y_n) = 1$ for $Y_n = 0$, and $\gamma(Y_n) = 0$ for $Y_n = 1$. The above equation gives rise to the well known bifurcation diagram with characteristic features such as period doubling cascades [e.g., *Feigenbaum*, 1978] for C larger than 3 and chaotic variations for C above ~3.57 (see, e.g., *Moon* [1987] and *Kuznetsov* [2004] for a detailed description). Such a behavior has been examined in a variety of contexts, from plasma physics [e.g., *Chen*, 1992] to biology [e.g., *Lesne*, 2006].

In the present report, instead of the above linear relationship between $\gamma(Y_n)$ and $Y_n$, we consider a nonlinear one and use the following functional form :

$$f(Y_n) = 0.5 \times \{1. - \tanh[A_Y(Y_n - 0.5)]\} \qquad (2)$$

Like in the linear case, this form $f(Y_n)$ varies between 1 for $Y_n = 0$ and 0 for $Y_n = 1$, the $A_Y$ coefficient controlling its maximum steepness. We then write the above growth rate γ as described in the introduction, viz.,

$$\gamma(Y_n) = \gamma_{min} + f(Y_n) \times [\gamma_{max} - \gamma_{min}] \qquad (3)$$

where $\gamma_{min}$ and $\gamma_{max}$ are the minimum and maximum values of $\gamma$, respectively. These limiting values are such that : $0 \leq \gamma_{min} < \gamma_{max} \leq 1$. Like in the logistic map (1), the recurrence relation now writes :

$$Y_{n+1} = C \times \gamma(Y_n) \times Y_n \qquad (4)$$

with C as a prescribed constant, and $\gamma(Y_n)$ as given in (3).

Figure 1 presents the growth rate $\gamma(Y_n)$ and the $Y_{n+1}$ profile for different values of $A_Y$, $\gamma_{min}$, and $\gamma_{max}$. In both panels of this figure, the profiles in dark blue and light blue are for moderate steepness of $f(Y_n)$ with $A_Y$ set to 5.0 in (2), and for two different minimum and maximum values of $\gamma$ (viz., 0.0 and 1.0 in dark blue, 0.1 and 0.9 in light blue). As for the red profile in Figure 1, it corresponds to a larger steepness of $f(Y_n)$ with $A_Y$ set to 7.5. These color coded profiles that exhibit smooth but nonmonotonous variations of both $\gamma$ and $Y_{n+1}$ as a function of $Y_n$, are to be compared with the dotted black lines in Figure 1 that show the linear case of the logistic map. The difference between linear and nonlinear growth rates is readily noticeable, which has an impact on the long-term behavior. In the following, this long-term behavior will be investigated via repeated iterations of the above recurrence relation. In particular, it will be shown that, unlike in the logistic map, the minimum value of $\gamma(Y_n)$ when $Y_n$ is maximum plays a prominent role in the structure of the bifurcation pattern.

3. Modeling Results

To examine the long term behavior of the above recurrence relation, we iterate this relation 1000 times to eliminate transients. We then plot the following 50 values of $Y_n$ as a function of coefficient C. We first consider cases with moderate steepness ($A_Y$ = 5.0 ; see blue profiles in Figure 1) and different values of $\gamma_{min}$ (the other parameter $\gamma_{max}$ being set to 1). The results obtained are shown in Figure 2. The top panel of this figure with $\gamma_{min}$ = 0 resembles the bifurcation diagram obtained for the logistic map as described by *May* [1976]. As mentioned above, period doubling cascades and transition to chaos clearly are noticeable in this figure as C increases from unity. Still, a notable difference between the pattern achieved in the top panel of Figure 2 and the well known bifurcation diagram is the overall shift toward smaller C values, with onset of period doubling near C = 1.75 and transition to chaos near C = 2.15.

In the center and bottom panels of Figure 2, the parameter $\gamma_{min}$ is set to 0.018 and 0.02, respectively. It can be seen in these panels that, although these minimum values of γ are very small, they considerably alter the structure of the bifurcation patterns. That is, instead of wide chaotic regions, narrow islands of confined chaos emerge for C above ~2.6, and the Y parameter only takes a limited number of values when C further increases toward 4. An example of these asymptotic $Y_n$ values is given in Figure 3 for C = 3.5 and for the three different values of $\gamma_{min}$ (0.0, 0.018, and 0.02, these values being coded in black, blue, and red, respectively). It is apparent from Figure 3 that period 3 is achieved for $\gamma_{min}$ = 0.02 (red profile), period 6 is achieved for $\gamma_{min}$ = 0.018 (blue profile), while no clear behavior is obtained for $\gamma_{min}$ = 0.0 (black profile).

Figure 4 presents another example of bifurcation patterns, considering now different values of $\gamma_{max}$ ($A_Y$ is unchanged, and $\gamma_{min}$ is set to 0.018 like in the center panel of Figure 2). Note that the change in $\gamma_{max}$ (by about 10%) that we consider here is much larger than that considered for $\gamma_{min}$ (less than 1%) in Figure 2. Again, a prominent change of the bifurcation pattern can be seen as compared to the pattern in the top panel of Figure 2, with a rapid evolution toward a small number (3 or 6) of $Y_n$ values for the largest C coefficients. Note also that the maximum $Y_n$ values (black dotted lines in each panel of Figure 4) nearly are distributed as $Y_n = 0.3 \times \gamma_{max} \times C$. Figures 2 and 4 thus demonstrate that, for given $A_Y$, the structure of the bifurcation pattern depends upon both $\gamma_{min}$ and $\gamma_{max}$ and actually exhibits an enhanced sensitivity to $\gamma_{min}$.

**4. Discussion**

To investigate the dependence of the recurrence relation upon $A_Y$, Figure 5 shows the results obtained for $A_Y = 7.5$. The general features discussed in Figure 2 with $A_Y = 5.0$ still are noticeable here. That is, a bifurcation pattern that resembles that of the logistic map can be seen in the top panel notwithstanding a substantial shift toward smaller C values. On the other hand, as $\gamma_{min}$ increases (center and bottom panels of Figure 5), narrow islands of $Y_n$ values are rapidly substituted to the wide chaotic regions. These islands reflect transition to periodic regimes with small periods (e.g., period 4 for C comprised between 2.6 and 2.8).

To provide further insights into the evolution of the bifurcation pattern depending upon $\gamma_{min}$, systematic computations were performed over a range of $\gamma_{min}$ values (from 0.0 to 0.2). For each value of $\gamma_{min}$ and for each value of C, we examined the computed $Y_n$ distribution with respect to a grid of 100 pixels regularly spaced between 0 and 1. In this grid, the number of pixels that contain computed $Y_n$ values gives an estimate of the $Y_n$ spread for the corresponding values of $\gamma_{min}$ and C. As an example, if a single asymptotic value of $Y_n$ is obtained as is the case in the C < 1.7 interval, only one pixel of the above grid is filled, yielding a spread of 1% which is much smaller than that in chaotic regions where one expects a large number of pixels to be filled. The results of these systematic computations are presented in Figure 6, considering either the linear case of the logistic map (left panel) or the nonlinear case with $A_Y = 5$ (center panel) or $A_Y = 7.5$ (right panel). In all cases, $\gamma_{max}$ is set to 1.

Figure 6 displays striking differences between the $Y_n$ spreads achieved for linear and nonlinear growth rates. In the former case, no dependence upon $\gamma_{min}$ is noticeable, except for a very slight change in (color coded) intensity as $\gamma_{min}$ increases. In other words, regardless of $\gamma_{min}$, the same bifurcation diagram prevails, with period doubling for C above ~3 and chaos for C above ~3.57. This situation clearly contrasts with that obtained using a nonlinear growth rate (center and right panels of Figure 6). In this latter case, the $Y_n$ spread considerably depends upon both $\gamma_{min}$ and C. If we identify the largest $Y_n$ spreads (in red) as chaotic regions, it can be seen that these regions are widely developed when $\gamma_{min}$ nears 0. As $\gamma_{min}$ increases, an erosion of these chaotic regions occurs and regimes with reduced $Y_n$ spread gradually develop. In other words, the increase of $\gamma_{min}$ tends to damp chaotic behaviors. In the right panels of Figure 6, numerous small-scale features also are noticeable that reflect prominent structuring of the bifurcation pattern.

## 5. Conclusions

The computations performed show that a nonlinear growth rate in place of a linear one in the logistic map leads to a different structure of the bifurcation pattern with shift and damping of chaotic regions. The long term behavior in this modified logistic map is found to critically depend upon both maximum steepness of the growth rate as well as its minimum and maximum values at, respectively, maximum and minimum values of the parameter of interest. This contrasts with the linear case where these parameters are not relevant. A prominent reordering of the logistic map may thus be obtained where chaotic regions are replaced by periodic regimes with low period numbers.

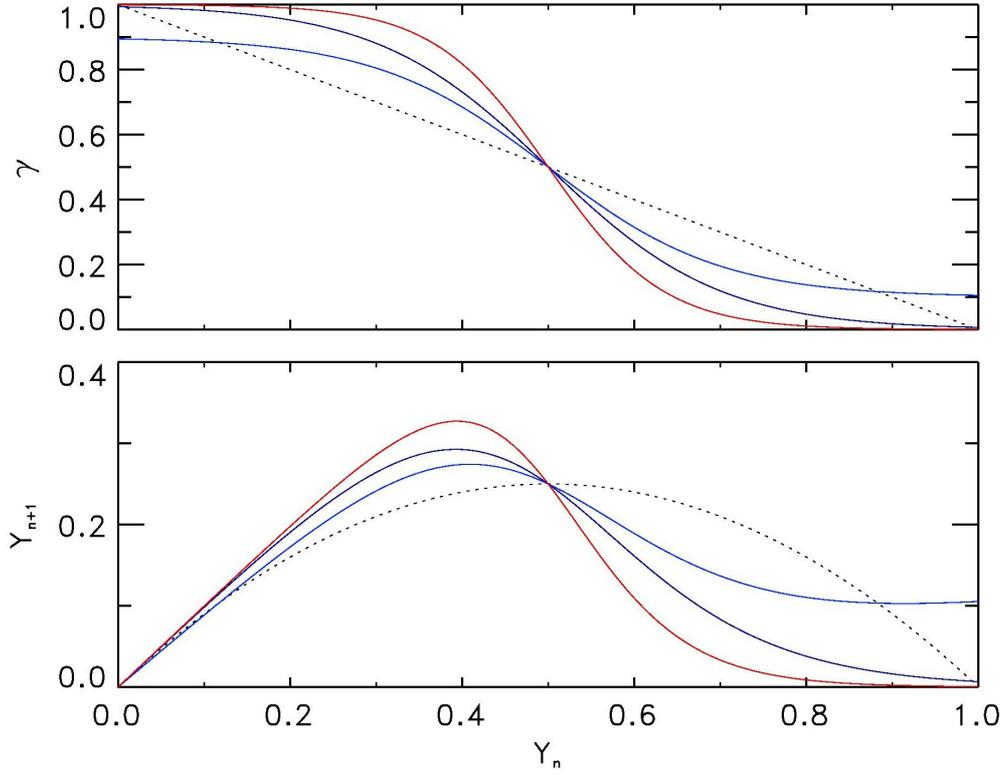

**Figure 1.** (*top*) Growth rate γ versus Y$_n$ and (*bottom*) Y$_{n+1}$ versus Y$_n$. Profiles in dark blue or light blue are for a hyperbolic tangent case with A$_Y$ = 5.0 and [γ$_{min}$, γ$_{max}$] = [0.0, 1.0] or [γ$_{min}$, γ$_{max}$] = [0.1, 0.9], respectively. The profile in red is for a hyperbolic tangent case with A$_Y$ = 7.5 and [γ$_{min}$, γ$_{max}$] = [0.0, 1.0]. Black dotted lines show the linear case of the logistic map.

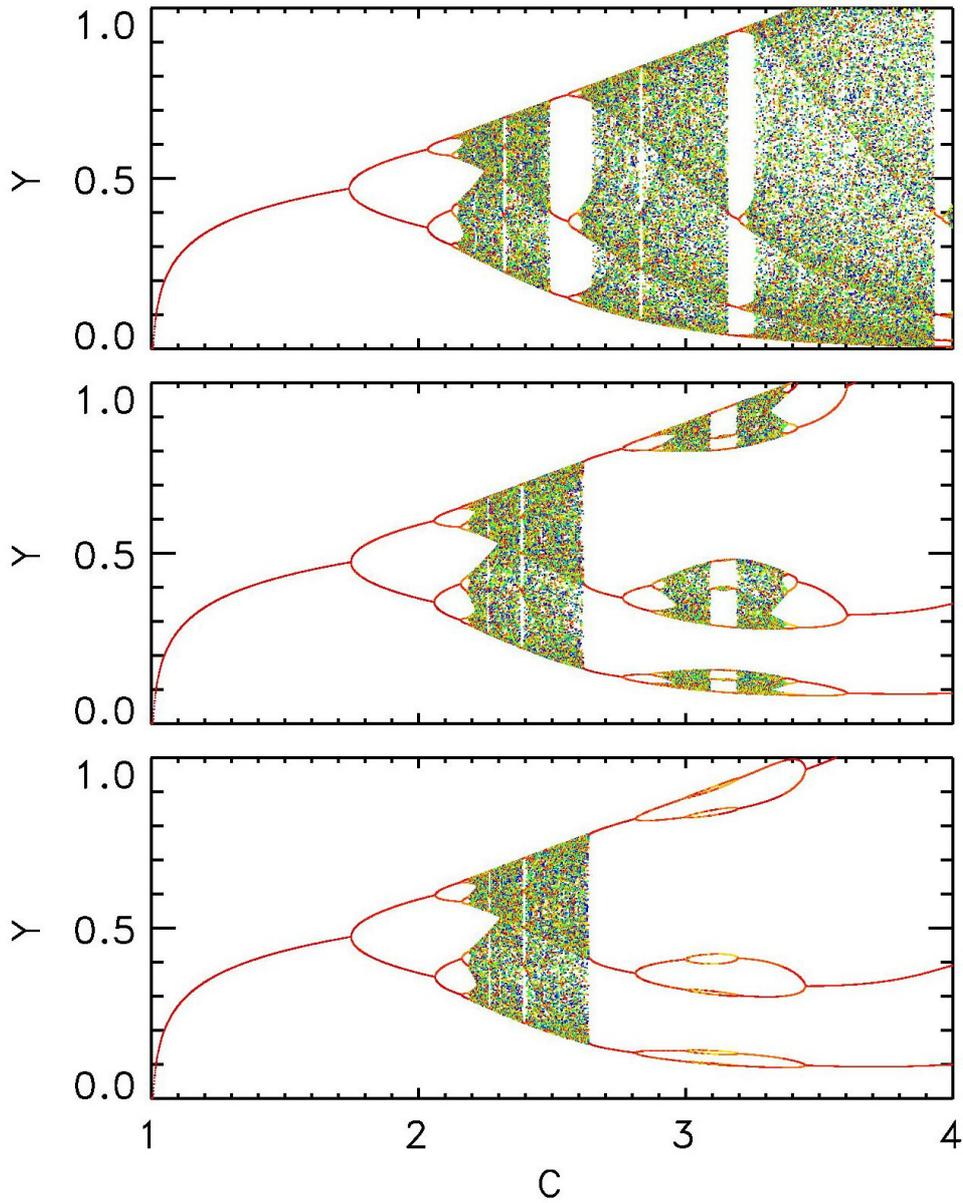

**Figure 2.** Bifurcation patterns (with iteration number color coded from blue to red) for a hyperbolic tangent case with $A_Y = 5.0$ and (*top*) $\gamma_{min} = 0.0$, (*center*) $\gamma_{min} = 0.018$, (*bottom*) $\gamma_{min} = 0.02$. In all three panels, $\gamma_{max}$ is set to 1.0.

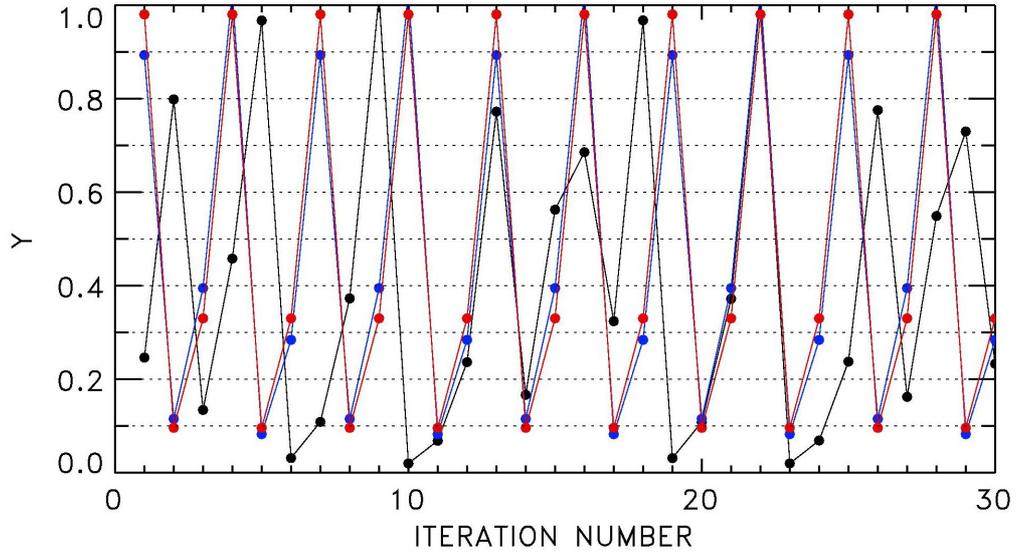

**Figure 3.** $Y_n$ versus iteration number for a hyperbolic tangent case with $A_Y = 5.0$ and $\gamma_{max} = 1$. The different colors correspond to $\gamma_{min} = 0.0$ (black), $\gamma_{min} = 0.018$ (blue), or $\gamma_{min} = 0.02$ (red). The coefficient C is set to 3.5.

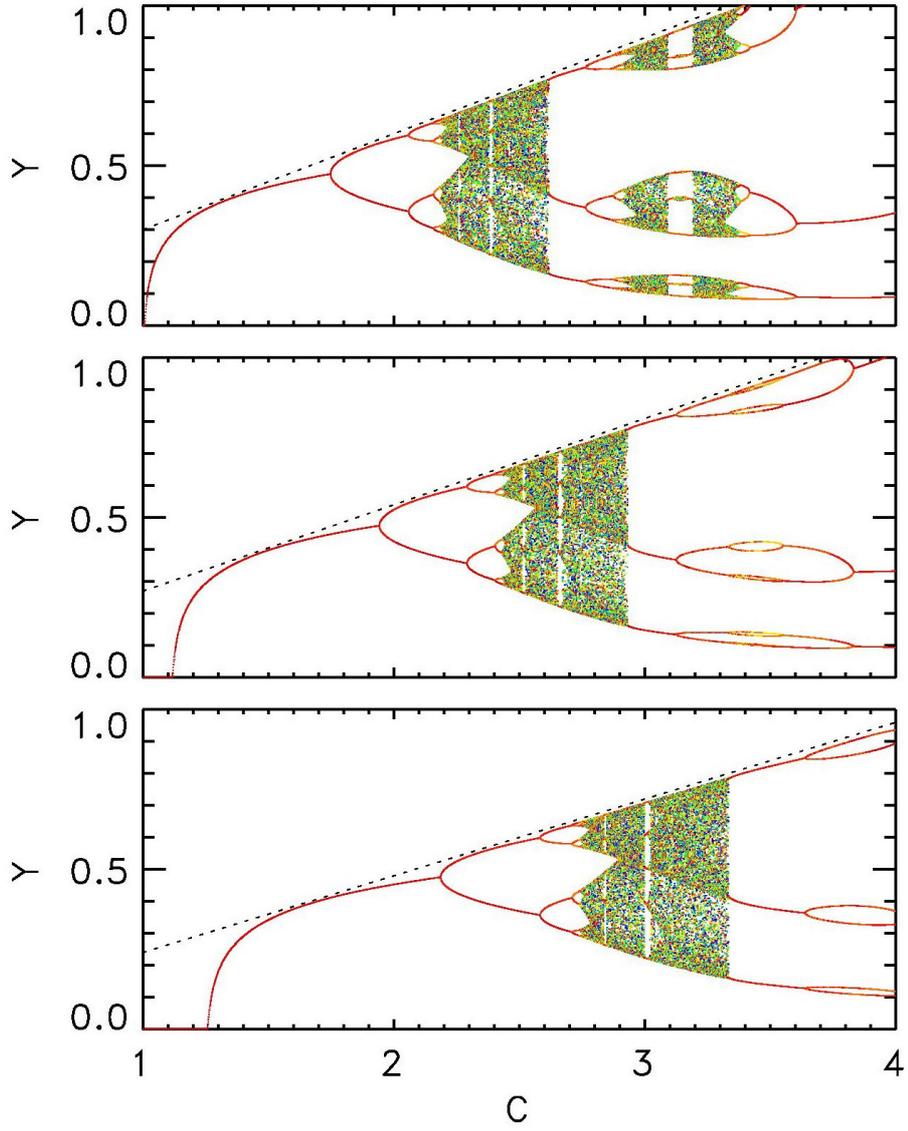

**Figure 4.** Identical to Figure 2 but with (*top*) $\gamma_{max} = 1.0$, (*center*) $\gamma_{max} = 0.9$, (*bottom*) $\gamma_{max} = 0.8$ in equation (3). Dotted lines show the $Y_n = 0.3 \times \gamma_{max} \times C$ profiles. In all three panels, $\gamma_{min}$ is set to 0.018.

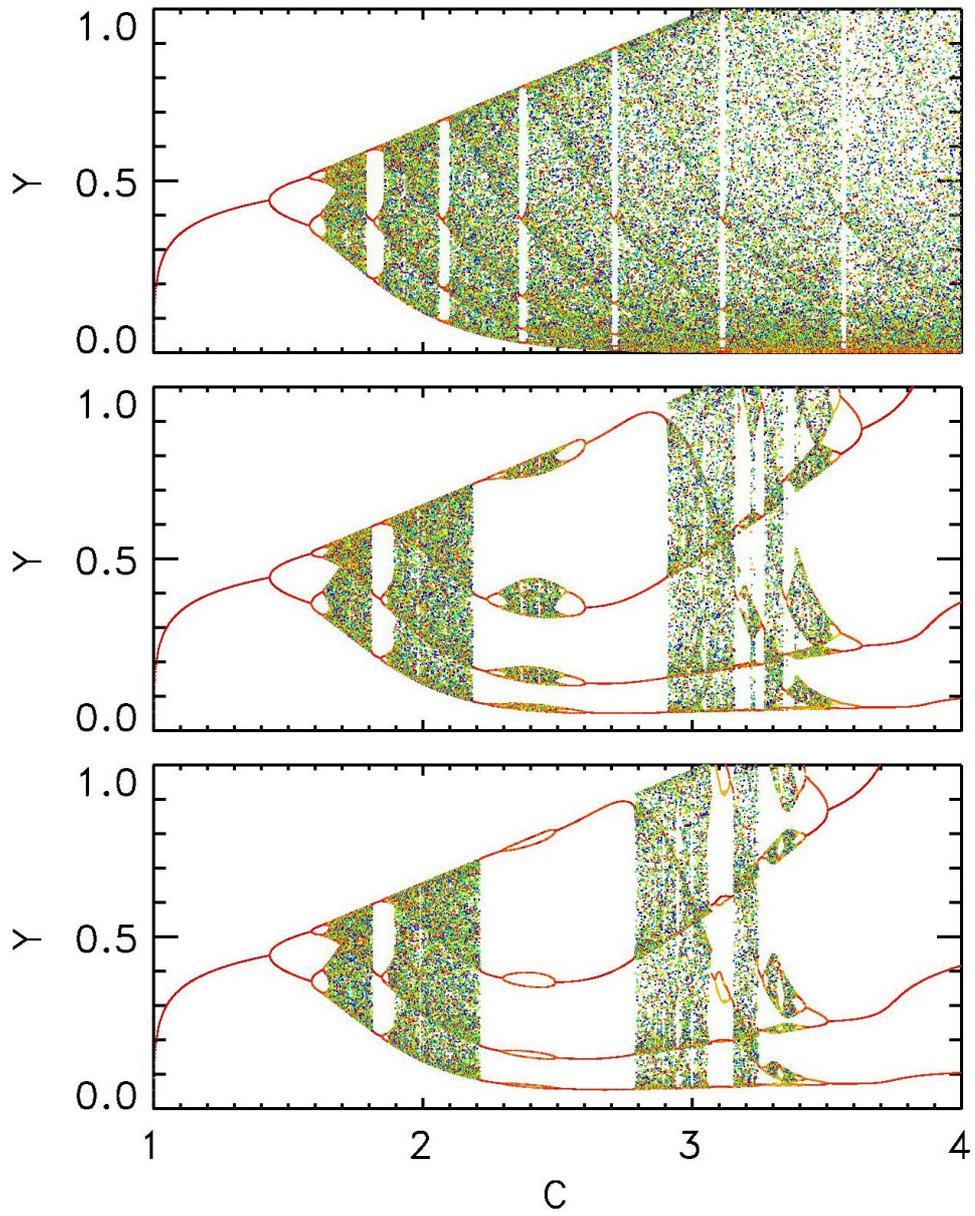

**Figure 5.** Identical to Figure 2 but with $A_Y$ set to 7.5 in equation (2).

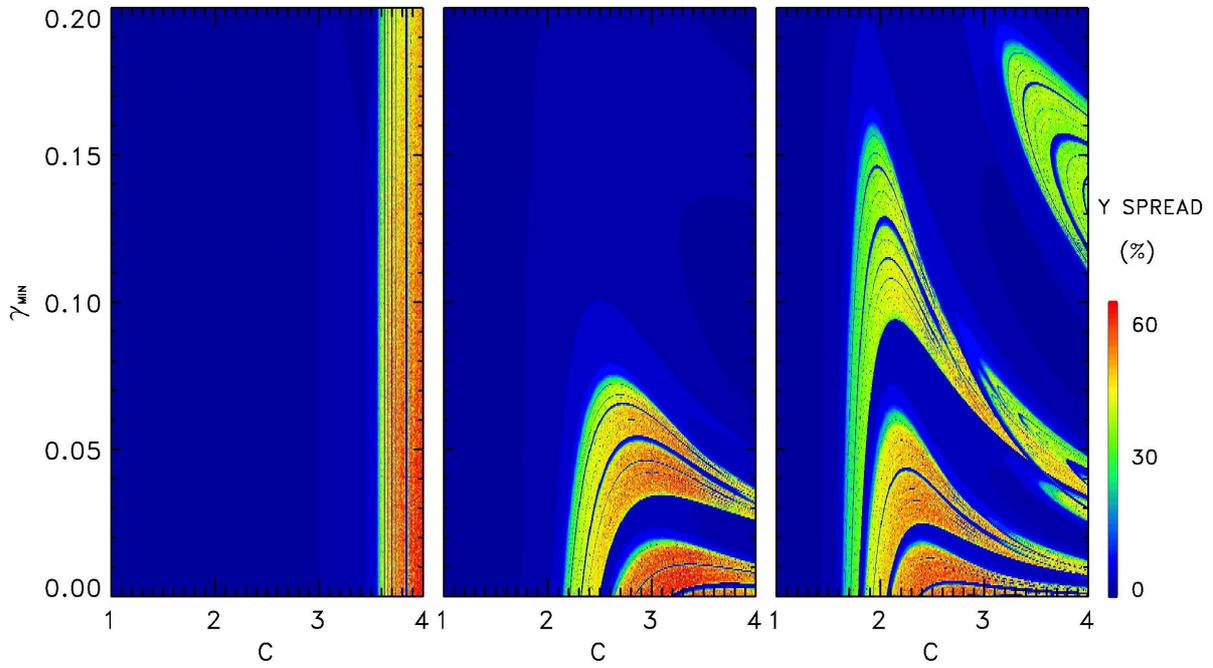

**Figure 6.** Color-coded $Y_n$ spread in the [0,1] interval as a function of C and $\gamma_{min}$. Three different cases are considered : (*left*) linear variation of $\gamma$ as in the logistic map, (*center and right*) hyperbolic tangent variation of $\gamma$ (equations (2)-(3)) with $A_Y = 5.0$ and $A_Y = 7.5$ in equation (2), respectively.